\documentclass[smallextended,envcountsect]{svjour3}

\let\endproof\relax

\usepackage{amsfonts}
\usepackage{amssymb}
\usepackage{amsthm}
\usepackage{amsmath}
\usepackage{stmaryrd}
\usepackage{float}
\usepackage{graphicx}
\usepackage{tikz}

\journalname{JOTA}
\def\citep{\cite}

\usetikzlibrary{shapes.geometric, fit, positioning, arrows}

\usepackage{hyperref}

\theoremstyle{plain}
\newtheorem{assumption}{Assumption}

\def\mathscr{\EuScript}

\newcommand{\EE}{\mathbb{E}}

\newcommand{\II}{\mathbb{I}}

\newcommand{\RR}{\mathbb{R}}

\newcommand{\UU}{\mathbb{U}}

\newcommand{\WW}{\mathbb{W}}
\newcommand{\XX}{\mathbb{X}}

\newcommand{\cE}{\mathcal{E}}

\newcommand{\cV}{\mathcal{V}}

{\normalsize}

{\normalsize}

{\normalsize}

\theoremstyle{plain}

\theoremstyle{definition}

\newcommand{\np}[1]{(#1)}                                   
\newcommand{\bp}[1]{\big(#1\big)}                           
\newcommand{\Bp}[1]{\Big(#1\Big)}                           
\newcommand{\bgp}[1]{\bigg(#1\bigg)}                        

\newcommand{\nc}[1]{[#1]}                                   
\newcommand{\bc}[1]{\big[#1\big]}                           
\newcommand{\Bc}[1]{\Big[#1\Big]}                           
\newcommand{\bgc}[1]{\bigg[#1\bigg]}                        

\newcommand{\na}[1]{\{#1\}}                                 
\newcommand{\ba}[1]{\big\{#1\big\}}                         

\newcommand{\projop}[1]{\mathrm{proj}_{#1}}                 

\makeatletter
\def\va@a{\boldsymbol{\va@arg^{\textstyle\text{\unboldmath$\scriptstyle\va@expo$}}_{\textstyle\text{\unboldmath$\scriptstyle\va@index$}}}}
\def\va#1{\def\va@expo{}\def\va@index{}\def\va@arg{\uppercase{#1}}%
  \@ifnextchar^{\va@h}{\@ifnextchar_\va@u\va@a}}
\def\va@h^#1{\def\va@expo{#1}\@ifnextchar_\va@hu\va@a}
\def\va@u_#1{\def\va@index{#1}\@ifnextchar^\va@uh\va@a}
\def\va@hu_#1{\def\va@index{#1}\va@a}
\def\va@uh^#1{\def\va@expo{#1}\va@a}
\makeatother

\newcommand{\feedback}{\gamma}                              

\newcommand{\price}{p}

\newcommand{\alloc}{r}

\newcommand{\horizon}{T}

\def\post{{t+1}}

\def\eqsepv{\; , \enspace}                                  
\def\eqfinv{\; ,}                                           
\def\eqfinp{\; .}                                           

\newcommand{\imag}{\mathrm{im}}

\newcommand{\argmin}{\mathop{\arg\min}}                     

\newcommand{\wrt}{\text{w.r.t.}}                            
\newcommand{\st}{\text{s.t.}}                               

\newcommand{\transp}{^{\top}}                               
\newcommand{\opt}{^{\sharp}}                                

\newcommand{\ic}[1]{\llbracket #1 \rrbracket}               

\makeatletter
\def\endproof{\finpreuve\@endtheorem}
\def\endremark{\finremark\@endtheorem}
\def\endexample{\finexemple\@endtheorem}
\makeatother

\newcommand{\finpreuvesymb}{$\Box$}                         
\newcommand{\finremarksymb}{$\Diamond$}                     
\newcommand{\finexemplesymb}{$\triangle$}                   

\newcommand{\finpreuve}{\ \hspace*{\fill}\finpreuvesymb}
\newcommand{\finremark}{\ \hspace*{\fill}\finremarksymb}
\newcommand{\finexemple}{\ \hspace*{\fill}\finexemplesymb}

\newcommand{\dynamic}{{g}}
\def\x{\boldsymbol{X}}
\def\u{\boldsymbol{U}}
\def\w{\boldsymbol{W}}

\def\Q{\va{Q}}
\def\F{\va{f}}

\def\pvaluefunc{\underline{V}}
\def\qvaluefunc{\overline{V}}

\def\adj{A}
\def\Adj{\mathcal{\adj}}

\newcommand{\kk}{{(k)}}
\newcommand{\kp}{{(k+1)}}

\def\horizon{T}
\def\final{T}

\newcommand{\CONE}{S}

\newcommand\produ{{\cV}}
\newcommand\transport{{\cE}}

\def\imag{\mathrm{im}}

\newcommand{\pscal}[2]{#1 \cdot #2}
\renewcommand{\pscal}[2]{\big\langle#1\:,#2\big\rangle}     
\def\tstep{\Delta T}
\newcommand{\VSDDP}{\underline{V}^{\mathrm{sddp}}}

\newcommand{\IGNORE}[1]{}                      

\newcommand{\FORALLTIMES}[3]{\forall #1\in\ic{#2,#3}}

\title{Upper and Lower Bounds for Large Scale Multistage Stochastic
       Optimization Problems: Application to Microgrid Management}

\titlerunning{Bounds for Large Scale
Multistage Stochastic Optimization Problems}

\author{P.~Carpentier \and J.-P.~Chancelier \and M.~De~Lara \and F.~Pacaud}

\institute{UMA, ENSTA Paris, IP Paris
\email{pierre.carpentier@ensta-paris.fr}
\and
Universit\'{e} Paris-Est, CERMICS (ENPC)
\email{chancelier@cermics.enpc.fr}
\and
Universit\'{e} Paris-Est, CERMICS (ENPC)
\email{delara@cermics.enpc.fr}
\and
\email{francois.pacaud@pm.me}
}

\date{\today}

\begin{document}

\setcounter{tocdepth}{2}
\setcounter{secnumdepth}{3}
\renewcommand{\labelitemi}{\textbullet}


\maketitle

\begin{abstract}
  We consider a microgrid where different prosumers exchange
  energy altogether by the edges of a given network.
  Each prosumer is located to a node of the network and
  encompasses energy consumption, energy production and storage
  capacities (battery, electrical hot water tank). The problem
  is coupled both in time and in space, so that a direct resolution
  of the problem for large microgrids is out of reach
  (curse of dimensionality). By affecting price or resources
  to each node in the network and resolving each nodal subproblem
  independently by Dynamic Programming, we provide decomposition
  algorithms that allow to compute a set of decomposed local value
  functions in a parallel manner. By summing the local value
  functions together, we are able, on the one hand, to obtain
  upper and lower bounds for the optimal value of the problem,
  and, on the other hand, to design global admissible policies
  for the original system. Numerical experiments are conducted
  on microgrids of different size, derived from data given
  by the research and development centre Efficacity,
  dedicated to urban energy transition.
  These experiments show that the decomposition algorithms give
  better results than the standard SDDP method, both in terms of
  bounds and policy values. Moreover, the decomposition methods are
  much faster than the SDDP method in terms of computation time,
  thus allowing to tackle problem instances incorporating more
  than 60 state variables in a Dynamic Programming framework.
\end{abstract}






\section{Introduction}

Multistage stochastic optimization problems are, by essence,
complex because their solutions are indexed both by stages
(time) and by uncertainties (scenarios). Hence, their large
scale nature makes decomposition methods appealing.
We refer to \citep{ruszczynski1997decomposition} and
\citep{carpentier2017decomposition} for a generic description
of decomposition methods in stochastic programming problems.
Dynamic Programming methods and their extensions are temporal
decomposition methods, that have been used on a wide panel of
problems, for example in dam management \citep{shapiro2012final}.
Spatial decomposition of large-scale optimization problems was
first studied in~\cite{cohen80}, and extended to open-loop
stochastic optimization problems \citep{cohenculioli90}.
Recent developments have mixed spatial decomposition methods with
Dynamic Programming to effectively solve large scale multistage
stochastic optimization problems. This work led to the introduction
of the so-called Dual Approximate Dynamic Programming (DADP)
algorithm, which was first applied to unit-commitment problems
with a single central coupling constraint linking different
stocks altogether~\citep{barty2010decomposition}.
We have extended this kind of methods in the companion
paper~\citep{bounds2019theory},
on the one hand by considering general coupling constraints
among units, and, on the other hand, by using two different
decomposition schemes, namely, price and resource decompositions.
This article presents applications of price and resource decomposition
schemes to the energy management of large scale urban microgrids.

General coupling constraints often arise from flows
conservation on a graph. Optimization problems on graphs
(monotropic optimization) have been studied since long
\citep{rockafellar1984network,bertsekas2008extended},
with applications, for example, to solve network utility
problems formulated as two-stage stochastic optimization
problems \citep{chatzipanagiotis2015augmented}.
Our motivation rather comes from electrical microgrid
management, where buildings (units) are able to consume,
produce and store energy and are interconnected through
a network. A broad overview of the emergence of consumer-centric
electricity markets is given in \cite{pinson2018emergence}.
We suppose here that all actors are benevolent, allowing
a central planner to coordinate the local units between
each other. Each local unit includes storages (hot water
tank and possibly a battery), and has to satisfy heat and
electrical demands. It also has the possibility to import
energy from an external regional grid if needed.
Some local units (prosumers) are able to produce their own
energy with solar panels, so as to satisfy their needs and
export the surplus to other consumers. The exchanges through
the network are modeled as a network flow problem on a graph.
We suppose that the system is impacted by uncertainties, both
in production (e.g. solar panels) or in demand (e.g. electrical
demands). Thus, the global problem can naturally be formulated
as a sum of local multistage stochastic optimization subproblems
coupled together via a global network constraint.
Such problems have been studied in \cite{thesepacaud}.
They are specially challenging from the dynamic optimization point
of view since the number of buildings may be large in a district.
We address districts with up to 48 buildings (with 64 associated
state variables), that is, a size largely beyond the limits
imposed by the well-known curse of dimensionality faced by
Dynamic Programming. The data associated with the districts
we are studying have been provided by Efficacity.
The local solar energy productions match realistic data
corresponding to a summer day in Paris. The local demands
are generated using a stochastic simulator experimentally
validated~\cite{schutz2015comparison}.
Efficacity is the urban Energy Transition Institute (ITE),
established in 2014 with the French government support.
The aim of Efficacity is to develop and implement innovative
solutions to build an energy-efficient and massively
carbon-efficient city.

\medskip

The paper is organized as follows.
In Sect.~\ref{chap:district:methods}, we model the global
optimization problem associated with a microgrid and apply
to it the main results obtained in the companion
paper~\cite{bounds2019theory}. We present both price and
resource decomposition schemes and recall how the Bellman
functions of the global problem are bounded above (resp. below)
by the sum of local resource-decomposed (resp. price decomposed)
value functions that satisfy recursive Dynamic Programming equations.
In Sect.~\ref{chap:district:numerics}, we present numerical
results for different microgrids of increasing size and complexity.
We compare the two decomposition algorithms with a state
of the art Stochastic Dual Dynamic Programming (SDDP) algorithm.
We analyse the convergence of all algorithms, and we compute the bounds
obtained by all algorithms. Thanks to the Bellman functions computed by
all algorithms, we are able to devise online policies for the initial
optimization problem and we compare the associated expected costs.
The analysis of case studies consisting of district microgrids
coupling up to 48 buildings together enlightens that decomposition
methods give better results in terms of economic performance, and
achieve up to a 4 times speedup in terms of computational time.


\section{Optimal management of a district microgrid}
\label{chap:district:methods}

In this section, we write the optimization problem corresponding
to a district microgrid energy management system on a graph
in \S\ref{sec:district:graphproblem}. We detail how to decompose
the problem node by node in \S\ref{sec:district:nodal},
both by using price and resource decomposition.
In \S\ref{sec:district:algorithm}, we show how to find the
most appropriate deterministic price and resource processes
for obtaining the best possible upper and lower bounds.

\subsection{Global optimization problem}
\label{sec:district:graphproblem}

A district microgrid is represented by a directed graph
$G = (\cV, \cE)$, with $\cV$ the set of nodes and $\cE$
the set of edges. We denote by $N_\cV$ the number of nodes,
and by $N_\cE$ the number of edges. Each node of the graph
corresponds to a single building comprising stocks, energy
production and local consumption. These buildings exchange
energy through the edges of the graph.

We first detail the different flows occurring in the graph
and the coupling constraints existing between flows in edges
and flows at nodes. We then formulate at each node a local
multistage stochastic optimization subproblem, as well as
a transportation subproblem on the graph. Finally, we gather
the coupling constraints and the subproblems inside a global
optimization problem.

\subsubsection{Exchanging flows through edges}

Flows are transported through the graph via the edges,
each edge~$e \in \ic{1, N_\cE}$ transporting a flow $q^e$
and each node $i \in \ic{1, N_\cV}$ importing or exporting
a flow $f^i$. Here $\ic{1, N}=\{1,2,\cdots,N\}$
denotes the set of integers between~1 and~$N$.

The node flows~$f^i$ and the edge flows~$q^e$
are related via a balance equation (Kirchhoff's current law),
which states that the sum of the algebraic edge flows arriving
at a particular node~$i$ is equal to the node flow $f^i$.
The Kirchhoff's current law can be written in matrix form as
\begin{equation}
  \label{eq:problem:nodearcflow}
  A q + f = 0 \eqfinv
\end{equation}
where $f = (f^1, \cdots, f^{N_\cV})\transp$ is the vector of
node flows, $q = (q^1, \cdots, q^{N_\cE})\transp$ is the vector
of edge flows and where $A \in \RR^{N_\cV \times N_\cE}$ is
the node-edge incidence matrix of the graph $G = (\cV, \cE)$.
Column~$e$ of~$A$ represents the edge~$e$ of the directed graph,
with values $+1$ (resp. $-1$) at the initial (resp. final) node
of the arc, and~$0$ elsewhere.

\subsubsection{Production cost on each node}

Each node of the graph~$G$ corresponds to a building which
may comprise stocks (hot water tank, battery), production
(solar panel), electric consumption. In case that the local
production cannot satisfy the local demand, external energy
is bought to the regional grid. We denote by $\final$ the time
horizon, by $\{0,1,\dots,\final\}$ the discrete time span (in
the application described in~\S\ref{Description_of_the_problems},
a unit period represents a 15mn time step). We write out all
random variables in bold. For a node~$i \in \ic{1, N_\cV}$,
the nodal subproblem is the minimization of a functional
$J_\produ^i\np{\F^i,x_0^i}$ depending on the node flow process
$\F^i = (\F_0^i, \cdots, \F^i_{\final-1})\transp$
arriving at node~$i$ between times $0$ and $\final-1$.

Let $\na{\XX_t^i}_{t\in \ic{0, \final}}$,
$\na{\UU_t^i}_{t\in \ic{0, \final-1}}$ and
$\na{\WW_t^i}_{t\in \ic{1, \final}}$
be sequences of Euclidian spaces of type~$\mathbb{R}^{p}$,
with appropriate dimensions~$p$ (possibly depending on time~$t$
and node~$i$). The optimal nodal cost $J_\produ^i$ is given by
\begin{subequations}
  \label{eq:district:localgenpb}
  \begin{align}
    J_\produ^i\np{\F^i,x_0^i} = \min_{\x^i, \u^i} \;
      & \EE \bgc{\sum_{t=0}^{\final-1}
                 L^i_t(\va X_t^i, \va U_t^i, \va W^i_{t+1}) +
                 K^i(\x_\horizon^i)}
        \label{eq:nodalcostequation} \eqfinv \\
    \st\ \FORALLTIMES{t}{0}{\final\!-\!1 }
      & \x_{t+1}^i = \dynamic_t^i(\x^i_t, \u_t^i, \w_{t+1}^i)
        \eqsepv \x_0^i = x_0^i
        \label{eq:dynamicsequation} \eqfinv \\
      & \Delta_t^i(\x_t^i, \u_t^i, \w_\post^i) = \F_t^i
        \label{eq:localbalanceequation} \eqfinv \\
      & \sigma(\u_t^i) \subset \sigma(\w_1, \cdots, \w_t, \w_\post)
        \label{eq:hazarddecisionframework} \eqfinv
  \end{align}
\end{subequations}
where we denote by $\x^i = (\x^i_0,\cdots,\x^i_\final)$,
$\u^i= (\u^i_0,\cdots,\u^i_{\final-1}) $
and $\w^i = (\w^i_0,\cdots,\w^i_\final)$ the local state
(stocks), control (production) and uncertainty (consumption)
processes. Constraint~\eqref{eq:localbalanceequation} represents
the energy balance inside node~$i$ for each time~$t$, with
$\Delta_t^i:\XX_t^i\times\UU_t^i\times\WW_\post^i\rightarrow\RR$.
In order to be able
to almost surely satisfy Constraints~\eqref{eq:localbalanceequation},
we assume that all decisions follow the \emph{hazard-decision}
information structure, that is, decision $\u_t^i$ is taken
after noise~$\w_{t+1}$ has been observed, hence the specific
form of Constraint~\eqref{eq:hazarddecisionframework}.
This slightly differs from the scope presented in the companion
paper~\cite{bounds2019theory} where the \emph{decision-hazard}
information structure was considered, but does not change
the kind of results already obtained.

We detail the dynamics~\eqref{eq:dynamicsequation} in building~$i$.
A battery is modeled with the linear dynamics
\begin{subequations}
  \begin{equation}
    \va b_\post^i = \alpha_b \va b_t^i + \tstep
    \bp{ \rho_c (\va u^{b,i}_t)^+ - \dfrac{1}{\rho_d} (\va u^{b, i}_t)^- }
    \eqsepv \forall t \in \ic{0,\dots, \final-1} \eqfinv
  \end{equation}
  where $\va b_t^i$ is the energy stored inside the battery
  at time $t$, $\u_t^{b,i}$ is the power exchanged with the battery,
  $\alpha_b$ is the auto-discharge rate and $(\rho_d, \rho_c)$
  are given yields. An electrical hot water tanks is modeled
  with the linear dynamics
  \begin{equation}
    \va h_\post^i = \alpha_h \va h_t^i + \tstep
    \bp{ \beta_h \u_t^{t,i} - \va d^{hw, i}_\post }
    \eqsepv \forall t \in \ic{0, .., \final-1} \eqfinv
  \end{equation}
  where $\va h_t^i$ is the energy stored inside the tank
  at time $t$, $\u_t^{t,i}$ is the power used to heat
  the tank, $\va d^{hw, i}_\post$ is the domestic hot water
  demand between time $t$ and $t+1$. The coefficient $\alpha_h$
  is a discharge rate corresponding to the losses by conduction
  and $\beta_h$ is a conversion coefficient.
\end{subequations}
Depending on the possible presence of a battery inside the building,
the nodal state $\x_t^i$ has dimension 2 or 1. If node $i$ has
a battery, its state is $\x_t^i = (\va b_t^i, \va h_t^i)$;
otherwise, its state is $\x_t^i = \va h_t^i$.
The value of the state at time $0$ is known, equal to $x_0^i$.

Equation~\eqref{eq:localbalanceequation} is the node balance
at node~$i$, with mapping~$\Delta_t^i$ given by
\begin{equation}
\Delta_t^i(\x_t^i, \u_t^i, \w_\post^i) =
\u_t^{ne,i} - \va{d}_\post^{el,i} - \u_t^{b,i} - \u_t^{t,i}
\eqfinv
\end{equation}
$\va d^{el, i}_\post$ being the residual\footnote{We have
chosen to aggregate the production of the solar panels of node~$i$
(if any) with the electricity demand, since they only appear by
their difference.} electricity demand between time $t$ and $t+1$,
and $\u_t^{ne,i}$ being the amount of electricity taken from
the external national grid.

Collecting the different variables involved in the model,
the control variable for building~$i$ is
$\u_t^i=(\u_t^{b,i},\u_t^{t,i},\u_t^{ne,i})$
and the noise variable affecting node~$i$ is
$\w_\post^i=(\va{d}_\post^{hw,i},\va{d}_\post^{el,i})$.

The cost function at node $i$ in \eqref{eq:nodalcostequation}
depends linearly on the price~$p_t^{el}$ to import electricity
from the external national grid, so that
\begin{equation}
L^i_t(\va X_t^i,\va U_t^i,\va W^i_{t+1}) = p_t^{el} \u_t^{ne,i}
\eqfinp
\end{equation}
The final cost~$K^i$ is a penalization to avoid an empty
electrical hot water tank at the end of the day.

The \emph{global nodal cost}
$J_\produ\bp{\F,\np{x_0^1,\ldots,x_0^{N_\cV}}}$
over the whole network is obtained by summing
the local nodal costs
\begin{equation}
  \label{eq:district:nodecriterion}
  J_\produ\bp{\F,\np{x_0^1,\ldots,x_0^{N_\cV}}}
= \sum_{i=1}^{N_\cV} J_\produ^i\np{\F^i,x_0^i} \eqfinp
\end{equation}

\subsubsection{Transportation cost on edges}

We now consider the edge costs arising when transporting
the flow~$\Q_t^e$ through each edge $e \in \ic{1, N_\cE}$
and for any time~$t\in\ic{0,\final-1}$.
The \emph{global edge cost}~$J_\transport(\Q)$ aggregates all
transport costs through the different edges in the graph, namely
\begin{equation}
  \label{eq:district:edgecriterion}
  J_\transport(\Q) =
  \EE \bgc{\sum_{e =1}^{N_\cE} \sum_{t=0}^{\final-1} l_t^e(\Q_t^e)}
  \eqfinv
\end{equation}
where $l_t^e: \RR \rightarrow \RR$ are convex real valued functions
assumed to be ``easy to compute'', e.g. quadratic. The cost $l_t^e$
can be induced by a difference in pricing, a fixed toll between
the different nodes, or by the energy losses through the network.

The global edge cost function $J_\transport$
in~\eqref{eq:district:edgecriterion}
is additive and thus decomposable \wrt\ (with respect to)
time and edges.

\subsubsection{Global optimization problem}

We have stated local nodal criteria~\eqref{eq:district:localgenpb}
and a global edge criterion~\eqref{eq:district:edgecriterion},
both depending on node and edge flows coupled by Constraint
\eqref{eq:problem:nodearcflow} at each time~$t\in \ic{0, \final-1}$,
that is, $\adj \Q_t + \F_t = 0$.
We rewrite these constraints in a single constraint involving
the global node flow and edge flow processes: $\Adj \Q + \F = 0$.
The matrix~$\Adj \in \RR^{\final . N_\cV \times \final . N_\cE}$
is a block-diagonal matrix with matrix~$\adj$ as diagonal element.
We are now able to formulate a \emph{global} optimization problem as
\begin{subequations}
  \label{transportproblem}
  \begin{align}
    V_0\opt(x_0^1,\cdots,x_0^{N_\cV}) = \min_{ \F,\Q} \;
      & \sum_{i=1}^{N_\cV} J_\produ^i\np{\F^i,x_0^i} + J_\transport(\Q) \\
    \st\
      & \Adj \Q + \F = 0
    \label{eq:couplingcons} \eqfinp
  \end{align}
\end{subequations}
Problem \eqref{transportproblem} couples $N_\cV + 1$ independent
criteria through Constraint~\eqref{eq:couplingcons}. As the resulting
criterion is additive and Constraint~\eqref{eq:couplingcons}
is affine, Problem~\eqref{transportproblem} has a nice form
to use decomposition-coordination methods.

\begin{remark}
There may be additional constraints in the problem,
for example bound constraints
$\underline F_t^i\leq\F_t^i\leq\overline F_t^i$
on the node flows, and bound constraints
$\underline Q_t^e\leq\Q_t^e\leq\overline Q_t^e$.
on the edge flows.
These constraints may be modeled,
in the global optimization problem criterion,
 by additional terms like
\begin{equation*}
\EE \bgc{\sum_{i=1}^{N_\cV}
         \II_{[\underline F_t^i,\overline F_t^i]^\final}(\F^i)} +
\EE \bgc{\sum_{e=1}^{N_\cE}\sum_{t=0}^{\final-1}
         \II_{[\underline Q_t^e,\overline Q_t^e]}(\Q_t^e)}
\eqfinv
\end{equation*}
where
\begin{equation*}
\II_{E}: x \mapsto
\left\{
\begin{aligned}
 0 \quad  & \text{ if } x \in E \\
 +\infty  & \text{ otherwise}
\end{aligned}
\right.
\end{equation*}
is the indicator function of the set~$E$. These additional terms
do not change the additive structure of the cost function.
\end{remark}

\subsection{Mixing nodal and time decomposition on a microgrid}
\label{sec:district:nodal}

In the companion paper~\cite{bounds2019theory}, we introduced
a generic framework to bound a global problem by decomposing
it into smaller local subproblems, easier to solve.
In Problem~\eqref{transportproblem}, the coupling
constraints~\eqref{eq:couplingcons} can be written
$(\F, \Q) \in - \CONE$, where the convex set
$\CONE$ of $\RR^{T \cdot N_\cV} \times \RR^{T \cdot N_\cE}$
is the linear subspace
\begin{equation}
  \label{eq:district:couplingcone}
  \CONE = \ba{(f, q) \in \RR^{T \cdot N_\cV} \times \RR^{T \cdot N_\cE}
  \; \big\vert \;
  \Adj q + f = 0 } \eqfinp
\end{equation}
Problem~\eqref{transportproblem} lies in the generic
framework introduced in~\cite{bounds2019theory}, and
the coupling equation $\Adj \Q + \F = 0$ becomes a special
case of the generic coupling constraint of this framework.
Moreover, it can easily be checked that the dual cone of
the set $S$ defined in~\eqref{eq:district:couplingcone}
has the following expression:
\begin{equation}
  \label{eq:conestarexpression}
  \CONE^\star = \ba{(\price, \mu) \in \RR^{T \cdot N_\cV} \times
  \RR^{T \cdot N_\cE}
  \; \big\vert \;
  \Adj^\top \price - \mu = 0 } \eqfinp
\end{equation}
The duality terms arising from Constraint~\eqref{eq:couplingcons}
are given by the formula
\begin{equation}
  \label{eq:district:convexconelambda}
  \pscal{\price}{f} + \pscal{\mu}{q} =
  \pscal{\price}{(\Adj q +f)} \eqsepv
  \forall (f, q) \in \CONE \eqsepv
  \forall (\price, \mu) \in \CONE^\star
  \eqfinv
\end{equation}
where $(u, v) \mapsto \pscal{u}{v}$ is the usual scalar
product on $\RR^{T \cdot N_\cV}$.
In order to solve Problem~\eqref{transportproblem}, we apply
the decomposition schemes introduced in \cite{bounds2019theory}.
More precisely, we first apply \emph{spatial decoupling}
into nodal and edge subproblems, and then apply the
\emph{temporal decomposition} induced by Dynamic Programming.

\subsubsection{Price decomposition of the global problem}

We follow the procedure introduced in \cite[\S2.2]{bounds2019theory}
to provide a lower bound and to solve Problem~\eqref{transportproblem}
by price decomposition.
We limit ourselves to deterministic price processes, that is, vectors
$\price=\np{\price^1, \ldots, \price^{ N_\cV}}\in\RR^{T \cdot N_\cV}$.
By Equation~\eqref{eq:district:convexconelambda},
the \emph{global price value function} associated with
Problem~\eqref{transportproblem} has the following expression, for all
$x_0=(x_0^1,\cdots,x_0^{N_\cV}) \in \XX_{0}^{1}\times\dots
 \times\XX_{0}^{N_\cV}$,
\begin{equation}
  \label{eq:district:dualvf}
  \pvaluefunc\nc{\price}\np{x_0} = \min_{\F, \Q} \;
  \sum_{i=1}^{N_\cV} J_\produ^i\np{\F^i,x_0^i} + J_\transport(\Q)
  + \EE \bc{\pscal{\price}{\Adj \Q + \F}} \eqfinp
\end{equation}
\begin{subequations}
\label{eq:district:localpricepb}
The global price value function~$\pvaluefunc\nc{\price}$ naturally
decomposes into a sequence of \emph{nodal price value functions}
\begin{equation}
  \label{eq:district:localnodalpricepb}
  \pvaluefunc_\produ^i\nc{\price^i}\np{x_0^i} = \min_{\F^i}
  J_\produ^i\np{\F^i,x_0^i} + \EE \bc{ \pscal{\price^i}{\F^i}}
  \eqsepv \forall i \in \ic{1, N_\cV} \eqfinv
\end{equation}
and an \emph{edge price value function} (which, to the difference
of the nodal price value function~\eqref{eq:district:localnodalpricepb},
does not depend on the initial state~$x_0$)
\begin{equation}
  \label{eq:district:localedgepricepb}
  \pvaluefunc_\transport\nc{\price} = \min_{\Q}\;
  J_\transport(\Q) + \EE \bc{ \pscal{\Adj^{\top}\price}{\Q}} \eqfinp
\end{equation}
\end{subequations}
For all $i \in \ic{1, N_\cV}$, considering the expression~\eqref{eq:district:localgenpb}
of the nodal cost~$J_\produ^i$, the nodal price value
function~\eqref{eq:district:localnodalpricepb} is,
for $x_0^i \in \XX_0^i$,
\begin{align*}
  \pvaluefunc_\produ^i\nc{\price^i}\np{x_0^i} =
  \min_{\x^i, \u^i, \F^i} \;
    & \EE \bgc{\sum_{t=0}^{\final-1}
      \Bp{L^i_t(\va X_t^i, \va U_t^i, \va W^i_{t+1}) +
      \pscal{\price_t^i}{\F_t^i}} +
      K^i(\x_\horizon^i)} \eqfinv \\
  \st
    & \:,\; \FORALLTIMES{t}{0}{\final\!-\!1} \eqsepv \nonumber \\
    & \x_{t+1}^i = \dynamic_t^i(\x^i_t, \u_t^i, \w_{t+1}^i)
      \eqsepv \x_0^i = x_0^i \eqfinv \\
    & \Delta_t^i(\x_t^i, \u_t^i, \w_\post^i) = \F_t^i \eqfinv \\
    & \sigma(\u_t^i) \subset \sigma(\w_1, \cdots, \w_t, \w_\post) \eqfinv
\end{align*}
The optimal value~$\pvaluefunc_\produ^i\nc{\price^i}(x_0^i)$
can be computed by Dynamic Programming under the so-called
white noise assumption.
\begin{assumption}
\label{hyp:independent}
The global uncertainty process $\np{\w_1,\cdots,\w_{\final}}$
consists of stagewise independent random variables.
\end{assumption}

For all node~$i \in \ic{1, N_\cV}$ and price $\price^i \in\RR^{T}$,
we introduce the sequence
$\na{\pvaluefunc_{\produ,t}^i\nc{\price^i}}_{t=0, \cdots, \final}$
of local price value functions defined,
for all $t \in \ic{0,\final}$ and $x_t^i \in \XX_t^i$, by
\begin{subequations}
\begin{align}
  \pvaluefunc_{\produ, t}^i\nc{\price^i}(x_t^i) =
  \min_{\x^i, \u^i, \F^i} \;
    & \EE \bgc{\sum_{s=t}^{\final-1}
      \Bp{L^i_s(\va X_s^i, \va U_s^i, \va W^i_{s+1}) +
      \pscal{\price_s^i}{\F_s^i}} +
      K^i(\x_\horizon^i)} \eqfinv \\
  \st
    & \:,\; \FORALLTIMES{t}{0}{\final\!-\!1} \eqsepv \nonumber \\
    & \x_{s+1}^i = \dynamic_s^i(\x^i_s, \u_s^i, \w_{s+1}^i)
      \eqsepv \x_t^i = x_t^i \eqfinv \\
    & \Delta_s^i(\x_s^i, \u_s^i, \w_{s+1}^i) = \F_s^i \eqfinv \\
    & \sigma(\u_s^i) \subset \sigma(\w_\post,\cdots,\w_s,\w_{s+1}) \eqfinv
\end{align}
\label{eq:bellmanpricefunctions}
\end{subequations}
with the convention $\pvaluefunc_{\produ, \final}^i\nc{\price^i} = K^i$.
Under Assumption~\ref{hyp:independent}, these local price
value functions satisfy the Dynamic Programming equations
for all $i \in \ic{1, N_\cV}$:
\begin{subequations}
\begin{align}
\pvaluefunc_{\produ, \final}^i\nc{\price^i}(x_\final^i)
 & = K^i(x_\final^i) \eqfinv \\
\intertext{and, for \( t=\final\!-\!1, \ldots, 0 \),}
\pvaluefunc_{\produ, t}^i\nc{\price^i}(x_t^i)
 & =  \EE\Bc{\min_{u_t^i} L^i_s(x_t^i,u_t^i, \w^i_{t})  +
             \pscal{\price_t^i}{\Delta_t^i(x_t^i,u_t^i,\w_{t+1}^i)}
             \nonumber \\
 & \hspace{3.3cm}
          + \pvaluefunc_{\produ, \post}^i\nc{\price^i}
            \bp{\dynamic_s^i(x_t^i,u_t^i, \w_{t+1}^i)}} \eqfinp
\end{align}
\label{eq:bellmanequationpricefunctions}
\end{subequations}
Note that the measurability constraints
$\sigma(\u_t^i) \subset \sigma(\w_1,\cdots,\w_{t+1})$ in the
above problem~\eqref{eq:bellmanpricefunctions} can be replaced
by $\sigma(\u_t^i) \subset \sigma(\w_1^i,\cdots,\w_\post^i)$
without changing the value
$\pvaluefunc_\produ^i\nc{\price^i}(x_0^i)$.
Indeed, Equation~\eqref{eq:bellmanpricefunctions} only involves
the local noise process $(\w_1^i,\cdots,\w_\final^i)$, so that
there is no loss of optimality to restrain the measurability
of the control process $\u^i$ to the filtration generated
by the local noise process $\w^i$.

Considering the expression~\eqref{eq:district:edgecriterion}
of the edge cost~$J_\transport(\Q)$,
the edge price value function~$\pvaluefunc_\transport\nc{\price}$
is additive \wrt\ time and space, and thus can be decomposed
at each time~$t$ and each edge~$e$. The resulting edge
subproblems do not involve any time coupling and can be computed
by standard mathematical programming tools or even analytically.

\subsubsection{Resource decomposition of the global problem}

We now solve Problem~\eqref{transportproblem} by resource
decomposition (see \cite[\S2.2]{bounds2019theory}) using
a deterministic resource process
$\alloc=\np{\alloc^1,\ldots,\alloc^{N_\cV}} \in \RR^{T \cdot N_\cV}$,
such that $\alloc \in \imag(\Adj)$.\footnote{If
$\alloc\notin\imag(\Adj)$, we have~$\qvaluefunc\nc{\alloc}=+\infty$ in~\eqref{eq:district:quantdec}
as the constraint $\Adj \Q + \alloc = 0$ cannot be satisfied.}
We decompose the global constraint~\eqref{eq:couplingcons}
\wrt\ nodes and edges as
\begin{equation*}
\F = \alloc \quad , \quad \Adj \Q = -\alloc \eqfinp
\end{equation*}
The \emph{global resource value function} associated to
Problem~\eqref{transportproblem} has the following expression,
for all
$x_0=(x_0^1,\cdots,x_0^{N_\cV}) \in
 \XX_{0}^{1}\times\dots\times\XX_{0}^{N_\cV}$,
\begin{subequations}
  \label{eq:district:quantdec}
  \begin{align}
    \qvaluefunc\nc{\alloc}\np{x_0} = \min_{\F, \Q} \;
      & \sum_{i=1}^{N_\cV} J_\produ^i\np{\F^i,x_0^i} + J_\transport(\Q) \\
    \st\
      & \F - \alloc = 0 \eqsepv \Adj \Q + \alloc = 0 \eqfinp
  \end{align}
\end{subequations}
\begin{subequations}
\label{eq:district:localallocpb}
The global resource value function~$\qvaluefunc\nc{\alloc}$ naturally
decomposes in a sequence of \emph{nodal resource value functions}
\begin{equation}
  \label{eq:district:localnodalallocpb}
  \qvaluefunc_\produ^i\nc{\alloc^i}\np{x_0^i} = \min_{\F^i} \; J_\produ^i\np{\F^i,x_0^i}
  \quad \st \quad \F^i - \alloc^i = 0
   \eqsepv \forall i \in \ic{1, N_\cV} \eqfinv
\end{equation}
and an \emph{edge resource value function}
(which does not depend on~$x_0$)
\begin{equation}
  \label{eq:district:localedgeallocpb}
  \qvaluefunc_\transport\nc{\alloc} = \min_{\Q} \; J_\transport(\Q)
  \quad \st \quad \Adj \Q + \alloc = 0 \eqfinp
\end{equation}
\end{subequations}
For all $i \in \ic{1, N_\cV}$, considering
the expression~\eqref{eq:district:localgenpb} of
the nodal cost~$J_\produ^i$, the nodal resource value
function~\eqref{eq:district:localnodalallocpb} is,
for $x_0^i \in \XX_0^i$,
\begin{align*}
  \qvaluefunc_\produ^i\nc{\alloc^i}(x_0^i) = \min_{\x^i, \u^i, \F^i} \;
    &  \EE \bgc{\sum_{t=0}^{\final-1}
                L^i_t(\va X_t^i, \va U_t^i, \va W^i_{t+1}) +
                K^i(\x_\horizon^i)} \eqfinv \\
  \st
    & \:,\; \FORALLTIMES{t}{0}{\final\!-\!1} \eqsepv \nonumber \\
    & \x_{t+1}^i = \dynamic_t^i(\x^i_t, \u_t^i, \w_{t+1}^i)
      \eqsepv \x_0^i = x_0^i \eqfinv \\
    & \Delta_t^i(\x_t^i, \u_t^i, \w_\post^i) = \F_t^i \eqfinv \\
    & \sigma(\u_t^i) \subset \sigma(\w_1, \cdots, \w_t, \w_\post) \eqfinv \\
    & \F^i_t - \alloc_t^i = 0 \eqfinp
\end{align*}
If Assumption~\ref{hyp:independent} holds true,
$\qvaluefunc_\produ^i\nc{\alloc^i}(x_0^i)$
can be computed by Dynamic Programming. That leads to
a sequence $\na{\qvaluefunc_{\produ,t}^i}_{t=0, \cdots, \final}$
of local resource value functions given, for all
$t \in \ic{0,\final}$ and $x_t^i \in \XX_t^i$, by
\begin{align*}
  \qvaluefunc_{\produ,t}^i\nc{\alloc^i}(x_t^i) = \min_{\x^i, \u^i, \F^i} \;
    &  \EE \bgc{\sum_{s=t}^{\final-1}
                L^i_s(\va X_s^i, \va U_s^i, \va W^i_{s+1}) +
                K^i(\x_\horizon^i)} \eqfinv \\
  \st
    & \:,\; \FORALLTIMES{s}{t}{\final\!-\!1} \eqsepv \nonumber \\
    & \x_{s+1}^i = \dynamic_s^i(\x^i_s, \u_s^i, \w_{s+1}^i)
      \eqsepv \x_t^i = x_t^i \eqfinv \\
    & \Delta_s^i(\x_s^i, \u_s^i, \w_{s+1}^i) = \F_s^i \eqfinv \\
    & \sigma(\u_s^i) \subset \sigma(\w_\post, \cdots, \w_s, \w_{s+1}) \eqfinv \\
    & \F^i_s - \alloc_s^i = 0 \eqfinv
\end{align*}
with the convention $\qvaluefunc_{\produ, \final}^i\nc{\alloc^i} = K^i$.
As already noticed in the case of price functions,
the measurability constraints
$\sigma(\u_s^i) \subset \sigma(\w_\post,\cdots,\w_{s+1})$
in the above problem can be replaced by the more restrictive
constraint $\sigma(\u_t^i) \subset \sigma(\w_1^i,\cdots,\w_\post^i)$
without changing the value $\qvaluefunc_\produ^i\nc{\alloc^i}(x_0^i)$.

In the case of resource decomposition, edges are coupled
through the constraint $\Adj \Q + \alloc = 0$, so that
the edge resource value function $\qvaluefunc_\transport\nc{\alloc}$
in \eqref{eq:district:localedgeallocpb} is not additive
in space, but remain additive \wrt\ time. As in price
decomposition, it can be computed by standard mathematical
programming tools or even analytically.

\subsubsection{Upper and lower bounds of the global problem}

Applying \cite[Proposition~2.2]{bounds2019theory}
to the global price value function~\eqref{eq:district:dualvf}
and resource value functions~\eqref{eq:district:quantdec},
we are able to bound up and down the optimal value $V_0\opt(x_0)$
of Problem~\eqref{transportproblem}, for all
$x_0=(x_0^1, \cdots, x_0^{N_\cV}) \in
 \XX_{0}^{1}\times\dots\times\XX_{0}^{N_\cV}$:
\begin{equation}
  \label{eq:district:upperlowerboundglobal}
  \sum_{i=1}^{N_\cV} \pvaluefunc_\produ^i\nc{\price^i}(x_0^i) +
               \pvaluefunc_\transport\nc{\price}
  \; \leq \; V_0\opt (x_0) \; \leq \;
  \sum_{i=1}^{N_\cV} \qvaluefunc_\produ^i\nc{\alloc^i}(x_0^i) +
  \qvaluefunc_\transport\nc{\alloc} \eqfinp
\end{equation}
From the expression \eqref{eq:conestarexpression} of the dual
cone $\CONE^{\star}$, which does not impose any constraint on
the vector $\price$, these inequalities hold true for any price
$\price \in \RR^{\horizon \cdot N_{\cV}}$, and for any resource
$\alloc \in \imag(\Adj)$.

\subsection{Algorithmic implementation}
\label{sec:district:algorithm}

In \S\ref{sec:district:nodal}, we decomposed
Problem~\eqref{transportproblem} spatially and temporally:
the global problem is split into (small) subproblems using
price and resource decompositions, and each subproblem is
solved by Dynamic Programming. These decompositions yield
bounds for the value of the global problem. To obtain
tighter bounds for the optimal value~\eqref{transportproblem},
we follow the approach presented in~\cite[\S3.2]{bounds2019theory},
that is, we maximize (resp. minimize) the left-hand side
(resp. the right-hand side)
in Equation~\eqref{eq:district:upperlowerboundglobal} \wrt\
the price vector~$\price\in\RR^{\horizon \cdot N_{\cV}}$
(resp. the resource vector~$\alloc\in\RR^{\horizon \cdot N_{\cV}}$).
We observe that determining optimal deterministic price and resource
coordination processes turns to implement gradient-like algorithms.

\subsubsection{Lower bound improvement}
\label{subsec:district:algorithm:price}

We detail how to improve the lower bound given
by the price value function
in~\eqref{eq:district:upperlowerboundglobal}. We fix
$x_0=(x_0^1,\cdots,x_0^{N_\cV}) \in
 \XX_{0}^{1}\times\dots\times\XX_{0}^{N_\cV}$, and we proceed
by maximizing the global price value function
$\pvaluefunc\nc{\price}\np{x_0}$ 
\wrt\ the deterministic price process~$\price$,
\begin{subequations}
\begin{equation}
  \label{eq:priceouterproblem}
  \sup_{\price\in \RR^{T \cdot N_\cV}} \;
  \pvaluefunc\nc{\price}\np{x_0} \eqfinv
\end{equation}
that is written equivalently (see Equation~\eqref{eq:district:dualvf})
\begin{equation}
  \label{eq:district:relaxedmpts}
  \sup_{\price} \; \min_{\F, \Q} \;
  \sum_{i=1}^{N_\cV} J_\produ^i\np{\F^i,x_0^i} + J_\transport(\Q)
  + \pscal{\price}{\EE\bc{\Adj \Q + \F}} \eqfinp
\end{equation}
\end{subequations}
We are able to maximize Problem~\eqref{eq:district:relaxedmpts}
\wrt\ $\price$ using a gradient ascent method (Uzawa algorithm).
At iteration~$k$, we suppose given a deterministic price process
$\price^\kk$ and a gradient step $\rho^\kk$.
The algorithm proceeds as follows:
\begin{subequations}
  \label{eq:district:subgradientpricedecom}
  \begin{align}
    \label{eq:flowupdate}
    {\F^i}^\kp
      & \in \argmin_{\F^i}
          J_\produ^i\np{\F^i,x_0^i} + \EE \bc{\pscal{{\price^i}^\kk}{\F^i}}
          \eqsepv \forall i \in \ic{1,N_\cV} \eqfinv \\
    \label{eq:edgeupdate}
    \Q^\kp
      & \in \argmin_{\Q}
          J_\transport(\Q) + \EE \bc{ \pscal{\Adj^{\top}\price^\kk}{\Q}}
          \eqfinv \\
    \price^\kp
      & = \price^\kk + \rho^\kk \; \EE\bc{\Adj \Q^\kp + \F^\kp} \eqfinp
    \label{eq:priceupdate}
  \end{align}
\end{subequations}
At each iteration $k$, updating $\price^\kk$ requires
the computation of the gradient
of~$\nabla \pvaluefunc\nc{\price^\kk}\np{x_0^1, \cdots, x_0^N}$,
that is, the expected value~$\EE\bc{\Adj \Q^\kp + \F^\kp}$,
usually estimated by Monte-Carlo.
The price update formula~\eqref{eq:priceupdate}, corresponding
to the standard gradient algorithm for the maximization
\wrt\ $\price$ in Problem~\eqref{eq:priceouterproblem},
can be replaced by more sophisticated methods (BFGS, interior
point method).

\subsubsection{Upper bound improvement}
\label{subsec:district:algorithm:alloc}

We now focus on the improvement of the upper bound given
by the global resource value function
in~\eqref{eq:district:upperlowerboundglobal}.
We fix $x_0=(x_0^1, \cdots, x_0^{N_\cV}) \in \XX_{0}^{1}\times\dots\times\XX_{0}^{N_\cV}$,
and we aim at solving the problem
\begin{subequations}
\begin{equation}
  \inf_{\alloc \in \imag(\Adj)}  \;
\qvaluefunc\nc{\alloc}\np{x_0} \; =
  \inf_{\alloc \in \imag(\Adj)} \sum_{i=1}^{N_\cV}
  \qvaluefunc_\produ^i\nc{\alloc^i}\np{x_0^i}
  + \qvaluefunc_\transport\nc{\alloc}
  \eqfinp
\end{equation}
The detailed expression of this problem is
(see Equation~\eqref{eq:district:quantdec})
\begin{multline}
  \label{eq:district:quantdecompositionrelaxed}
  \inf_{\alloc \in \imag(\Adj)} \bgp{ \sum_{i=1}^{N_\cV}
  \Bp{\min_{\F^i} \; J_\produ^i\np{\F^i,x_0^i}
  \quad \st \quad \F^i - \alloc^i \!=\! 0} \\
  + \Bp{\min_{\Q} \; J_\transport(\Q)
  \quad \st \quad \Adj \Q + \alloc \!=\! 0}} \eqfinp
\end{multline}
\end{subequations}

The gradients \wrt\ $\alloc$, namely
$\mu^i = \nabla_{\alloc^i} \qvaluefunc_\produ^i\nc{\alloc^i}\np{x_0^i}$
and $\xi = \nabla_{\alloc} \qvaluefunc_\transport\nc{\alloc}$,
are obtained when computing the nodal resource value
functions~\eqref{eq:district:localnodalallocpb} and the edge
resource value function~\eqref{eq:district:localedgeallocpb}.
The minimization problem~\eqref{eq:district:quantdecompositionrelaxed}
is then solved using a gradient-like method. At iteration $k$,
we suppose given the resource $\alloc^\kk$ and a gradient
step $\rho^\kk$. The algorithm proceeds as follows:
\begin{subequations}
  \label{eq:district:subgradientquantdecom}
  \begin{align}
    \label{eq:flowquantupdate}
    {\F^i}^\kp
      & \in \argmin_{\F^i}
        J_\produ^i\np{\F^i,x_0^i} \quad \st \quad \F^i - \alloc^i
        \eqsepv \forall i \in \ic{1,N_\cV} \eqfinv \\
    \label{eq:edgequantupdate}
    \Q^\kp
      & \in \argmin_{\Q}
        J_\transport(\Q) \quad \st \quad \Adj \Q + \alloc \!=\! 0
        \eqfinv \\
    \alloc^\kp
      & = \projop{\imag(\Adj)} \Bp{\alloc^\kk - \rho^\kk \;
        \bp{\mu^\kp + \xi^\kp}} \eqfinv
    \label{eq:district:quantupdate}
  \end{align}
where $\projop{\imag(\Adj)}$ is the orthogonal projection
onto the subspace $\imag(\Adj)$. Again, the projected gradient
algorithm \eqref{eq:district:quantupdate} used to update
the resource can be replaced by more sophisticated methods.
\end{subequations}

\section{Application to microgrids optimal management}
\label{chap:district:numerics}

In this section, we treat an application. We apply the price and resource
decomposition algorithms described in \S\ref{sec:district:algorithm}
to a microgrid management problem, where different buildings are
connected together. The energy management system (EMS) controls
the different energy flows inside the microgrid, so as to ensure
at each node and at each time that the production meets the demand
at least cost. We give numerical results comparing the price and
resource decomposition algorithms with the Stochastic Dual Dynamic
Programming (SDDP) algorithm.

\subsection{Description of the problems}
\label{Description_of_the_problems}

We look at a microgrid connecting different buildings together.
As explained in~\S\ref{sec:district:graphproblem}, we model
the distribution network as a directed graph with buildings
set on nodes and distribution lines set on edges. The buildings
exchange energy with each other via the distribution network.
If the local production is unable to fulfill the local demand,
energy can be imported from an external regional grid
as a recourse.

The network configuration corresponds to heterogeneous domestic
buildings. Each building is equipped with an electrical hot water
tank, some have solar panels and some others have batteries. As batteries
and solar panels are expensive, they are shared out across the network.
We view batteries and electrical hot water tanks as energy stocks.
Depending on the presence of battery inside the building, the state $\x_t^i$
at node~$i$ has dimension 2 or 1 (energy stored inside the water tank and
energy stored in the battery), and such is the control $\u_t^i$ at node~$i$
(power used to heat the tank and power exchanged with the battery).
Furthermore, we suppose that all agents are benevolent
and share the use of their devices across the network.

We limit ourselves to a one day horizon. We look at a given day
in summer, discretized at a 15mn time step, so that $\final = 96$.
Each house has its own electrical and domestic hot water
demand profiles. At node~$i$, the uncertainty~$\w_t^i$
is a two-dimensional vector, namely the local electricity demand
and the domestic hot water demand. We choose to aggregate
the production of the solar panel with the local electricity
demand. We model the distribution of the uncertainty~$\w_t^i$
with a finite probability distribution on the set~$\WW_t^i$.

\begin{figure}[!ht]
  \centering
  \begin{tabular}{ccc}
    \includegraphics[width=3cm]{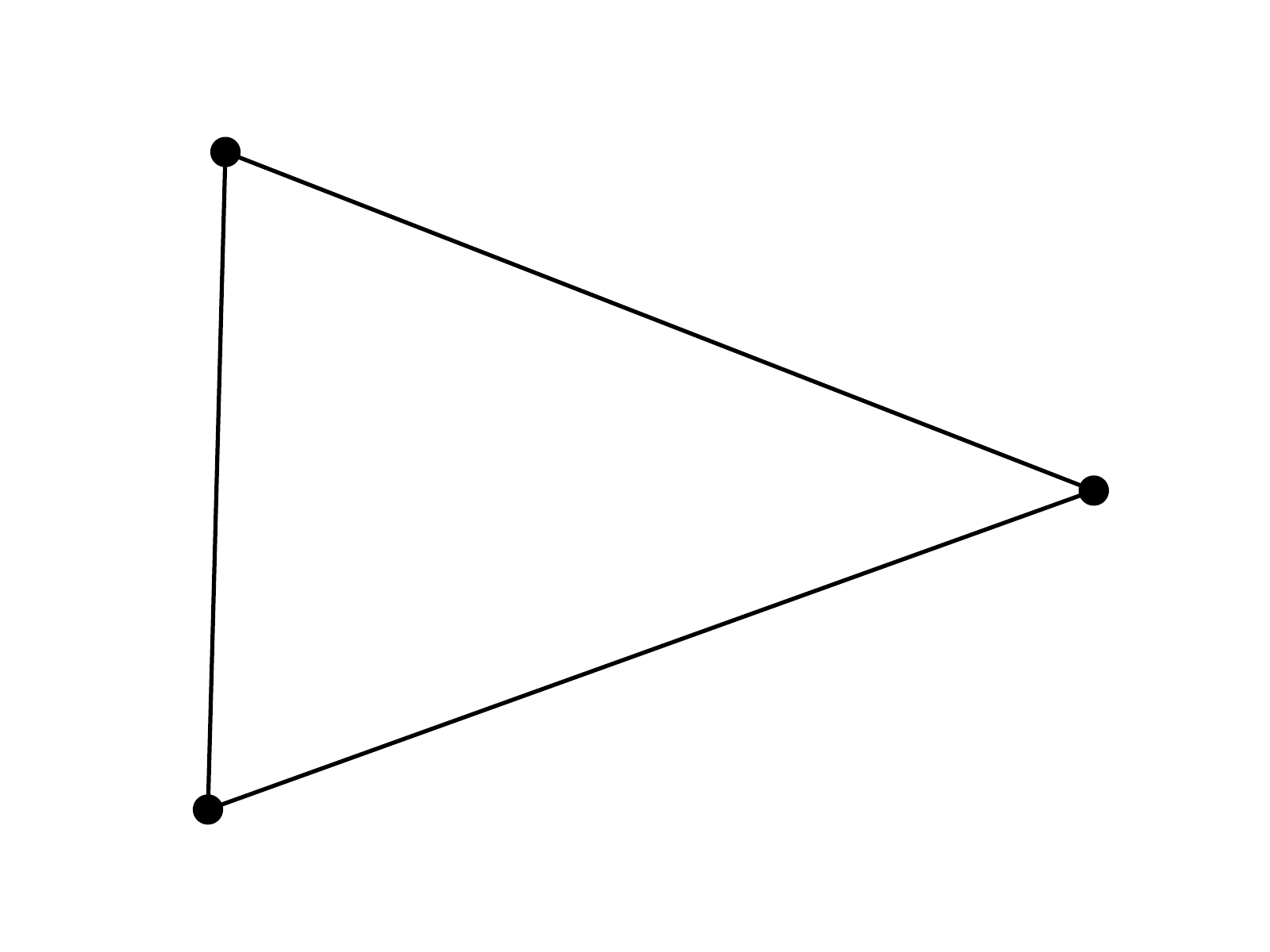}&
    \includegraphics[width=3cm]{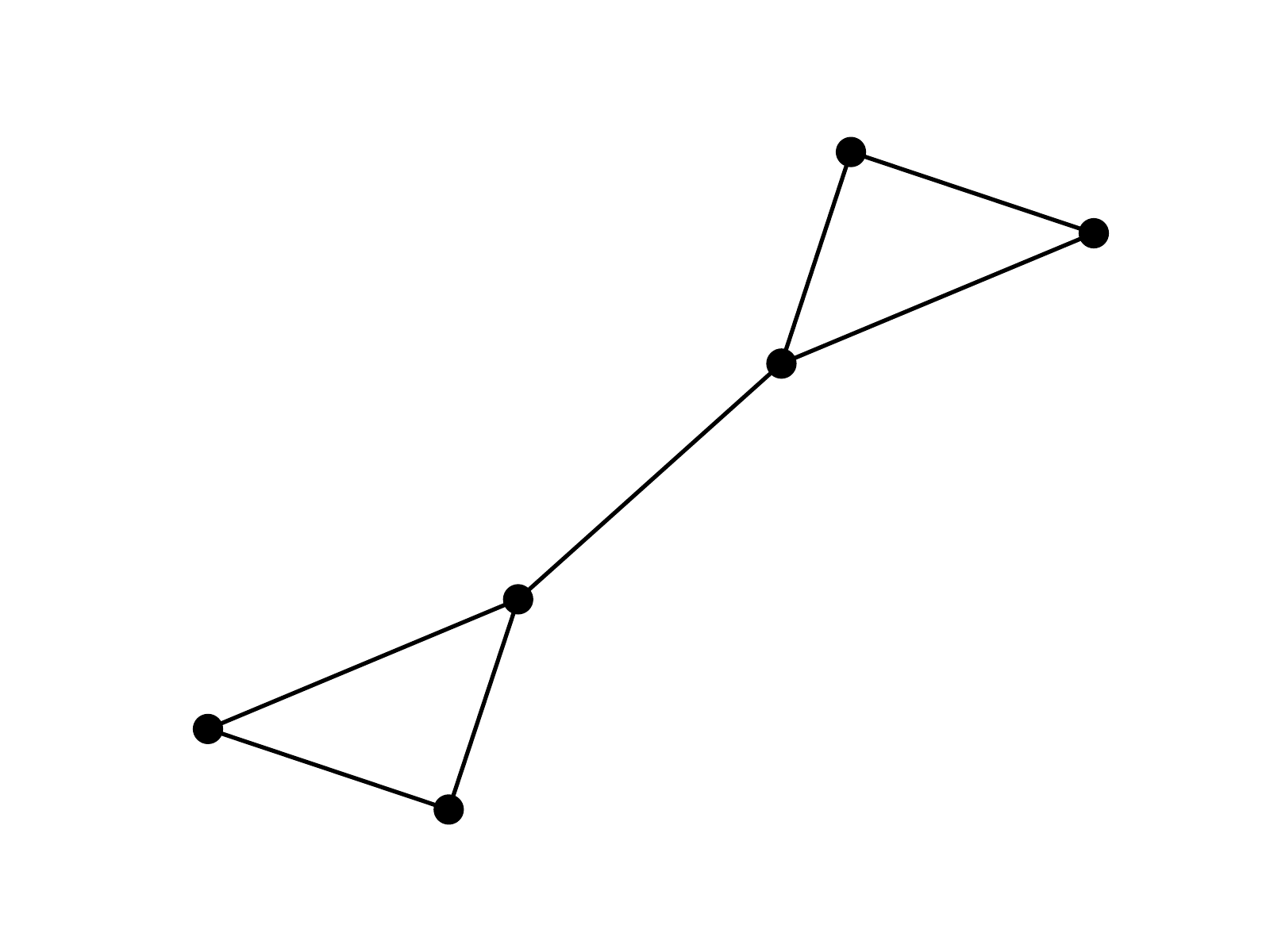} &
    \includegraphics[width=3cm]{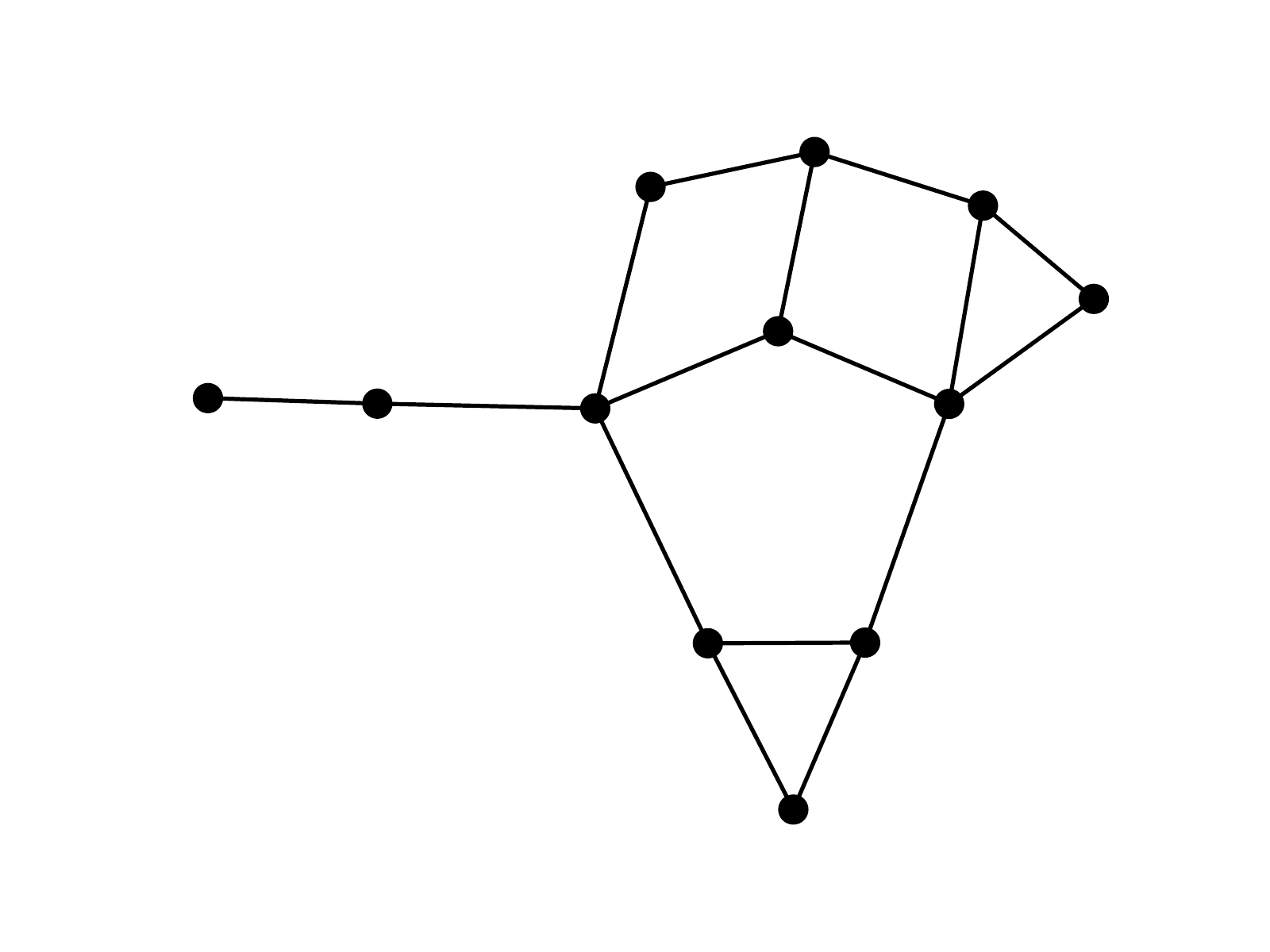} \\
    \textrm{3-Nodes} &
    \textrm{6-Nodes} &
    \textrm{12-Nodes}
  \end{tabular}
  \begin{tabular}{cc}
    \includegraphics[width=4.5cm]{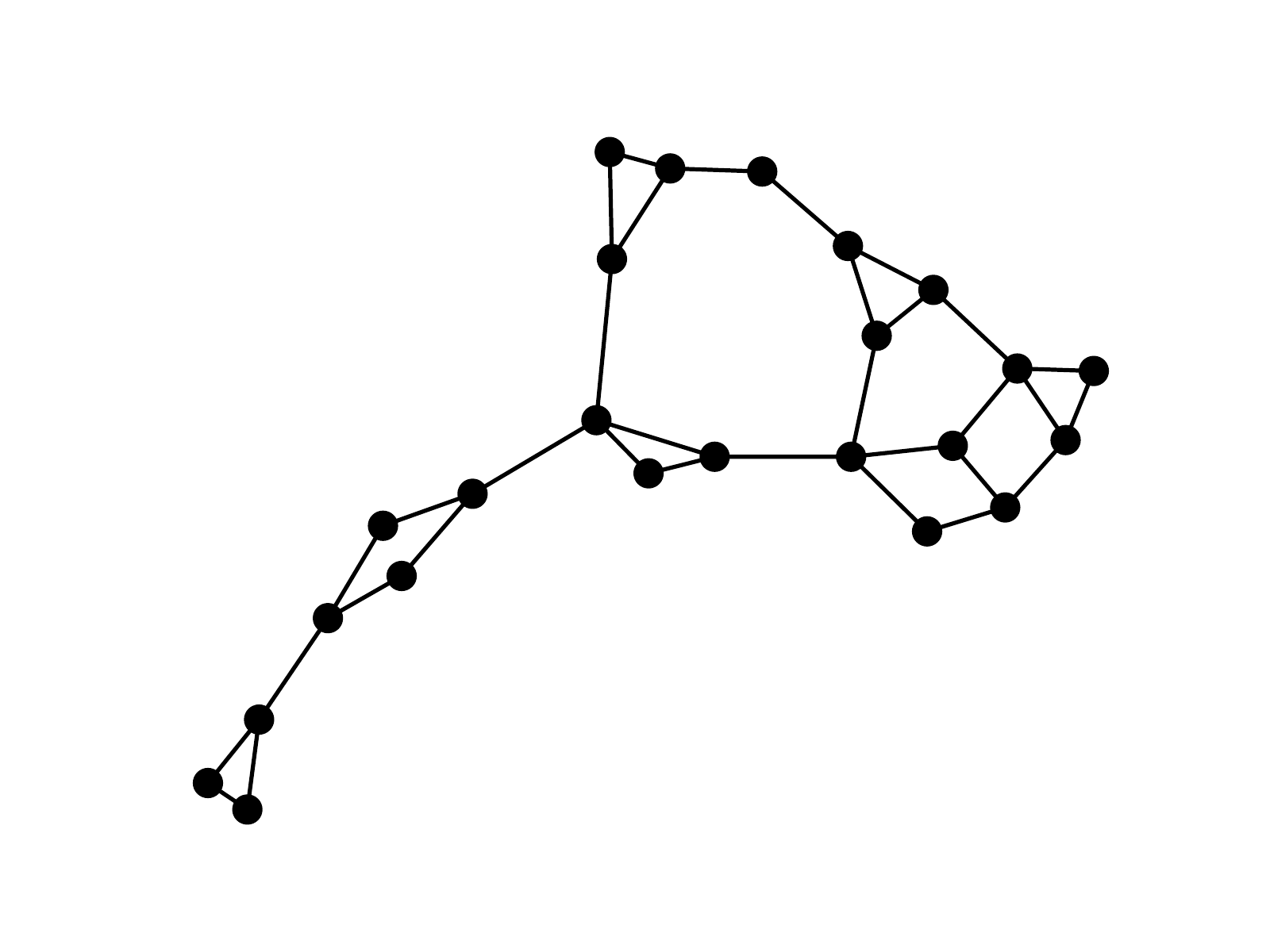} &
    \includegraphics[width=4.5cm]{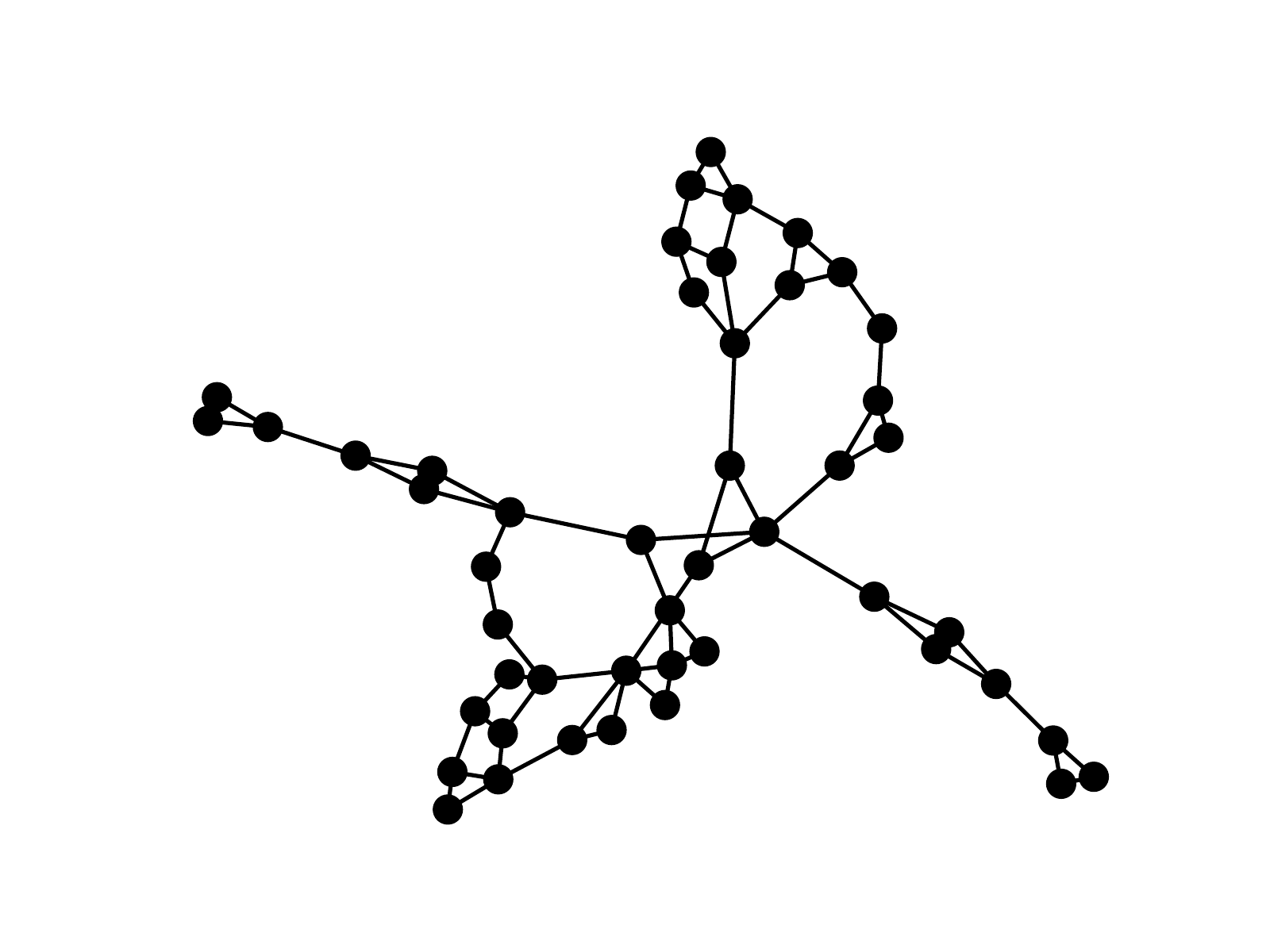} \\
    \textrm{24-Nodes} & \textrm{48-Nodes}
  \end{tabular}
  \caption{Topologies of the different graphs connecting
           buildings in the microgrids}
  \label{fig:nodal:graphtopology}
\end{figure}

We consider six different problems with growing sizes.
Table~\ref{tab:numeric:pbsize} displays the different dimensions
considered.
\begin{table}[!ht]
  \centering
  {\normalsize
  \begin{tabular}{|c|ccccc|}
    \hline
    Problem           & $N_\cV$ (nodes) & $N_\cE$ (edges) & $dim(\XX_t)$ & $dim(\WW_t)$ & $supp(\w_t)$ \\
    \hline
    \hline
    \textrm{3-Nodes}  & 3               & 3               & 4            & 6            & $10^3$    \\
    \textrm{6-Nodes}  & 6               & 7               & 8            & 12           & $10^6$    \\
    \textrm{12-Nodes} & 12              & 16              & 16           & 24           & $10^{12}$ \\
    \textrm{24-Nodes} & 24              & 33              & 32           & 48           & $10^{24}$ \\
    \textrm{48-Nodes} & 48              & 69              & 64           & 96           & $10^{48}$ \\
    \hline
  \end{tabular}
  }
  \caption{Microgrid management problems with growing dimensions}
  \label{tab:numeric:pbsize}
\end{table}
As an example, the 12-Nodes problem consists of twelve buildings;
four buildings are equipped with a 3~kWh battery, and four other
buildings are equipped with 16~$\text{m}^2$ of solar panels.
The devices are dispatched so that a building equipped with a solar
panel is connected to at least one building with a battery.
The support size of each local random variable~$\w_t^i$ remains low,
but that of the global uncertainty $\w_t=(\w_t^1, \cdots, \w_t^{N_\cV})$
becomes huge as~$N_\cV$ grows, so that the exact computation
of an expectation \wrt\ $\w_t$ is out of reach. The topologies of
the different graphs are depicted in Figure~\ref{fig:nodal:graphtopology}.
The structure of the microgrid as well as the repartition of batteries
and solar panel on it come from case studies provided by
the urban Energy Transition Institute Efficacity.

\subsection{Resolution algorithms}
\label{ssec:nodalalgorithms}

We reconsider the two decomposition
algorithms introduced in~\S\ref{sec:district:algorithm} and apply them
to each problem described in Figure~\ref{fig:nodal:graphtopology}.
We will term \emph{Dual Approximate Dynamic Programming} (DADP)
the price decomposition algorithm described
in \S\ref{subsec:district:algorithm:price} and
\emph{Primal Approximate Dynamic Programming} (PADP) the resource
decomposition algorithm described in \S\ref{subsec:district:algorithm:alloc}.
We compare DADP and PADP with the well-known Stochastic Dual Dynamic
Programming (SDDP) algorithm (see~\cite{girardeau2014convergence}
and references inside) applied to the global problem.

\subsubsection{Gradient-like algorithms}

It is common knowledge that
the usual gradient descent algorithm may be slow to converge.
To overcome this issue, we use a quasi-Newton algorithm to
approximate numerically the Hessian of the two global value functions
$\price \mapsto \pvaluefunc\nc{\price}\np{x_0^1,\cdots,x_0^{N_\cV}}$ in~\eqref{eq:district:dualvf} and
$\alloc \mapsto \qvaluefunc\nc{\alloc}\np{x_0^1, \cdots, x_0^{N_\cV}}$ in~\eqref{eq:district:quantdec}.
More precisely, the quasi-Newton algorithm is performed using
Ipopt 3.12 compiled with the MUMPS linear solver
(see~\citep{wachter2006implementation}). The algorithm stops either
when a stopping criterion is fulfilled or when no descent direction
is found.

\subsubsection{SDDP on the global problem}
\label{SDDP_on_the_global_problem}

In order to have at disposal a reference solution for the global
problem~\eqref{transportproblem}, we solve it using the Stochastic
Dual Dynamic Programming (SDDP) method. But the SDDP algorithm is
not implementable in a straightforward manner. Indeed,
the cardinality of the global noise support becomes huge with
the number of nodes~$N_\cV$ (see Table~\ref{tab:numeric:pbsize}),
so that the exact computation of expectations, as required at each time
step during the backward pass of the SDDP algorithm (see~\cite{shapiro10}),
becomes untractable.
To overcome this issue, we resample the probability distribution
of the global noise~$\np{\w_t^1,\cdots,\w_t^{N_\cV}}$ for each time~$t$
to deal with a noise support of reasonable size. To do so, we use the
$k$-means clustering method, as described in \cite{rujeerapaiboon2018scenario}.
By using the Jensen inequality \wrt\ the noises, we know that the optimal
quantization of a finite distribution yields a new optimization problem
whose optimal value is a lower bound for the optimal value of the original
problem, provided that the local problems are convex \wrt\ the noises
(see \cite{lohndorfmodeling} for details). Then, the exact lower bound
given by SDDP with resampling remains a lower bound for the exact lower
bound given by SDDP without resampling, which istself is a lower bound
for the original problem by construction.
In the numerical application, we fix the resampling size to $100$.
We denote by $\na{\VSDDP_t}_{t=0, \cdots, T}$
the value functions returned by the SDDP algorithm.
Notice that, whereas the SDDP algorithm suffers from
the cardinality of the global noise support, the DADP and PADP
algorithms do not.

We stop SDDP when the gap between its exact lower bound and a statistical
upper bound is lower than 1\%. That corresponds to the standard SDDP's
stopping criterion described in~\cite{shapiro10}, which is reputed
to be more consistent than the first stopping criterion introduced
in \cite{pereira1991multi}.
SDDP uses a level-one cut selection algorithm~\citep{guigues2017dual}
and keeps only the 100 most relevant cuts. By doing so, we significantly
reduce the computation time of SDDP.

\subsection{Devising control policies}
\label{ssec:nodalpolicies}

Each algorithm (DADP, PADP and SDDP) returns a sequence of value
functions indexed by time, that allow to build a global control
policy. Using these value functions, we define a sequence of
\emph{global value functions} $\na{\widehat{V}_t}_{t\in \ic{0, \final}}$ approximating the original value functions:
\begin{itemize}
  \item
    $\displaystyle \widehat{V}_t =
    \VSDDP_t$ \hspace{0.0cm} for SDDP,
  \item
    $\displaystyle \widehat{V}_t =
     \sum_{i=1}^{N_\cV} \pvaluefunc_{\produ,t}^i\nc{\price}
     + \pvaluefunc_{\transport, t}\nc{\price}$ for DADP,
  \item
    $\displaystyle \widehat{V}_t =
     \sum_{i=1}^{N_\cV} \qvaluefunc_{\produ,t}^i\nc{\alloc}
     + \qvaluefunc_{\transport, t}\nc{\alloc}$ for PADP.
\end{itemize}
We use these global value functions to build a global control policy
for all time $t\in \ic{0, \final-1}$. For any global state~$x_t \in \XX_t$
and global noise $w_\post \in \WW_\post$,
the control policy is a solution of the following one-step DP problem:
\begin{subequations}
  \label{eq:district:admissiblestrategy}
  \begin{align}
    \feedback_t(x_t, w_\post) \in \argmin_{u_t} \min_{f_t, q_t}
      & \sum_{i=1}^{N_\cV} L_t^i(x_t^i, u_t^i, w^i_\post) \! + \!
        \sum_{e=1}^{N_\cE} l_t^e(q_t^e) \! + \!
        \widehat{V}_\post\bp{x_\post}  \\
    \st\
      & x_\post^i = \dynamic_t^i(x_t^i, u_t^i, w^i_\post)
        \eqsepv \forall i\in \ic{1, N_\cV} \eqfinv \\
      & \Delta_t^i(x_t^i, u_t^i, w_\post^i) = f_t^i
        \eqsepv \forall i\in \ic{1, N_\cV} \eqfinv \\
      & A q_t + f_t = 0 \eqfinp
  \end{align}
\end{subequations}
As the strategy induced by~\eqref{eq:district:admissiblestrategy}
is admissible for the global problem~\eqref{transportproblem},
the expected value of its associated cost is an upper bound
of the optimal value~$V_0\opt$ of the original minimization
problem~\eqref{transportproblem}.

\subsection{Numerical results}
We first compare the three algorithms depicted in
\S\ref{ssec:nodalalgorithms}. We analyze the convergence
of them and the CPU time needed for achieving it.
We also present the value of the exact bounds obtained
by each algorithm. Then we evaluate the quality of
the strategies~\eqref{eq:district:admissiblestrategy}
introduced in \S\ref{ssec:nodalpolicies} for
the three algorithms.

\subsubsection{Computation of the Bellman value functions}

\begin{figure}[!ht]
  \centering
  \includegraphics[width=10.0cm]{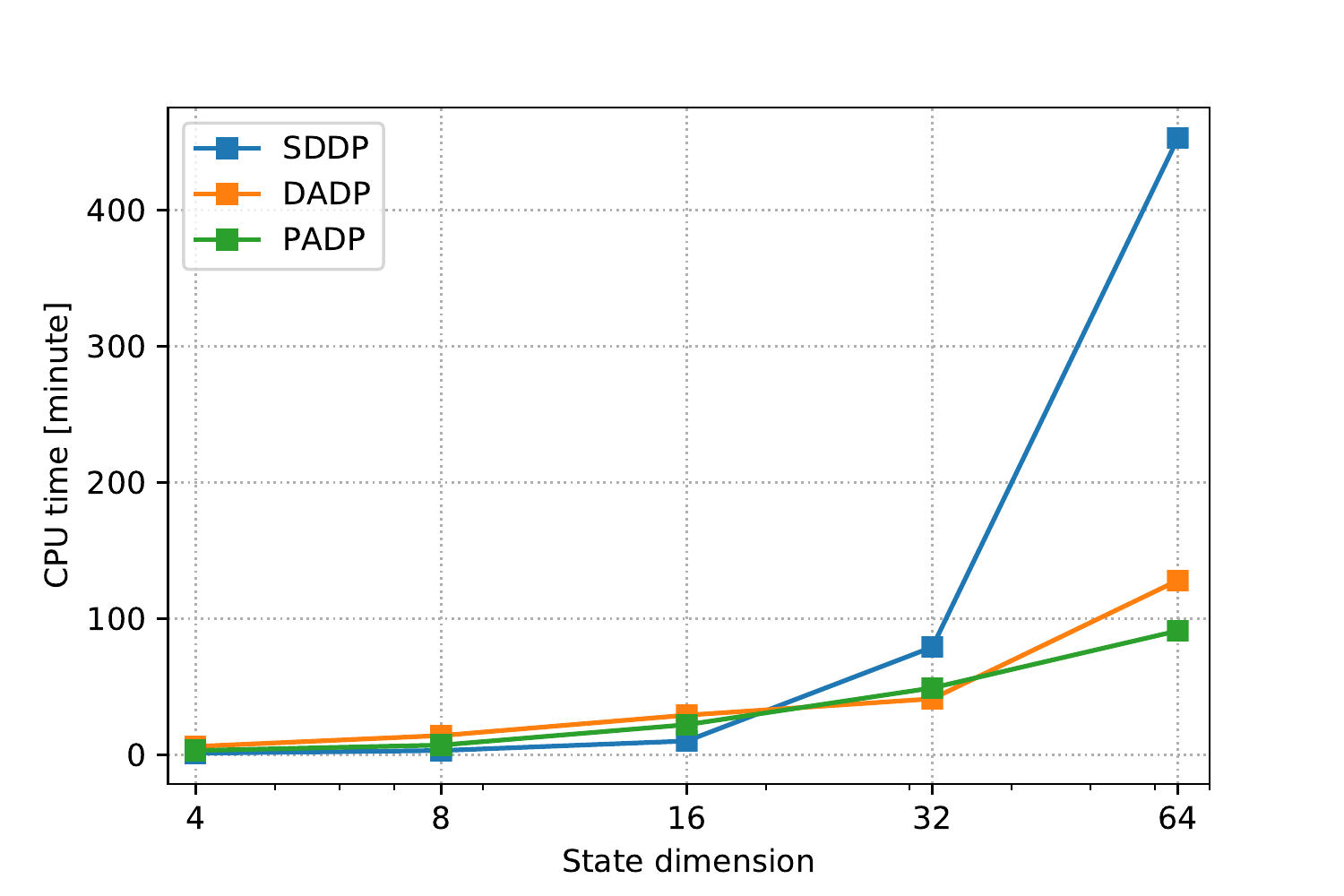} \\
  \caption{CPU time for the three algorithms as a function
           of the state dimension}
  \label{fig:nodal:cputime}
\end{figure}

We solve Problem~\eqref{transportproblem} by SDDP,
price decomposition (DADP) and resource decomposition (PADP).
Table~\ref{tab:district:numeric:optres} details the execution
time and number of iterations taken before reaching convergence.
\begin{table}[!ht]
  \centering
  {\normalsize
  \begin{tabular}{|l|ccccc|}
    \hline
    Problem            & \textrm{3-Nodes}  \hspace{-0.2cm}
                       & \textrm{6-Nodes}  \hspace{-0.2cm}
                       & \textrm{12-Nodes} \hspace{-0.2cm}
                       & \textrm{24-Nodes} \hspace{-0.2cm}
                       & \textrm{48-Nodes} \hspace{-0.2cm} \\
    \hline
    $|\XX_t|$          & 4                & 8                & 16
                       & 32               & 64     \\
    \hline
    \hline
    SDDP CPU time      & 1'               & 3'               & 10'
                       & 79'              & 453'   \\
    SDDP iterations    & 30               & 100              & 180
                       & 500              & 1500   \\
    \hline
    \hline
    DADP CPU time      & 6'               & 14'              & 29'
                       & 41'              & 128'   \\
    DADP iterations    & 27               & 34               & 30
                       & 19               & 29     \\
    \hline
    \hline
    PADP CPU time      & 3'               & 7'               & 22'
                       & 49'              & 91'    \\
    PADP iterations    & 11               & 12               & 20
                       & 19               & 20     \\
    \hline
  \end{tabular}
  }
  \caption{Convergence results for SDDP, DADP and PADP}
  \label{tab:district:numeric:optres}
\end{table}
For a small-scale problem like \textrm{3-Nodes} (second column
of Table~\ref{tab:district:numeric:optres}), SDDP is faster
than DADP and PADP. However, for the 48-Nodes problem (last
column of Table~\ref{tab:district:numeric:optres}),
\emph{DADP and PADP} are \emph{more than three times faster}
than SDDP. Figure~\ref{fig:nodal:cputime} depicts how much CPU
time take the different algorithms with respect to the number
of state variables of the district. For this case study, we observe
that the \emph{CPU time grows almost linearly} \wrt\ the number
of nodes for DADP and PADP, whereas it grows exponentially for SDDP.
Otherwise stated, decomposition methods scale better than SDDP
in terms of CPU time for large microgrids instances.

\begin{figure}[!ht]
  \centering
  \includegraphics[width=7cm]{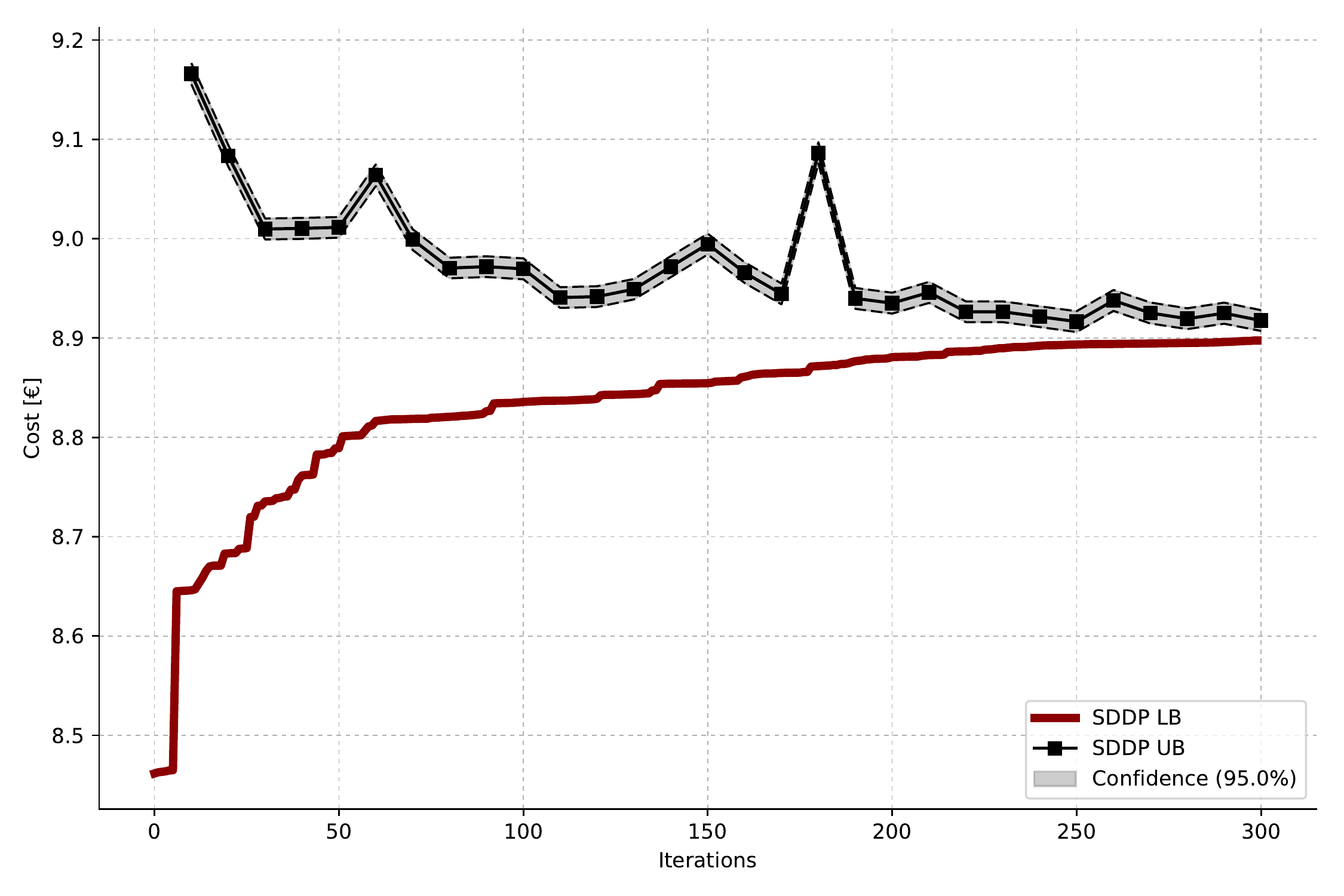}
  \caption{Evolution of SDDP lower and upper bounds for the 12-Nodes
  problem}
  \label{fig:numerics:convsddp}
\end{figure}

\paragraph{Convergence of the SDDP algorithm.}
Figure~\ref{fig:numerics:convsddp} displays the convergence of SDDP
for the 12 nodes problem. The approximate upper bound is estimated
every 10 iterations, with 1,000 scenarios. We observe that the gap
between the upper and lower bounds is below 1\% after 180 iterations.
The lower bound remains stable after 250 iterations.

\begin{figure}[!ht]
  \centering
  \begin{tabular}{cc}
    \includegraphics[width=5cm]{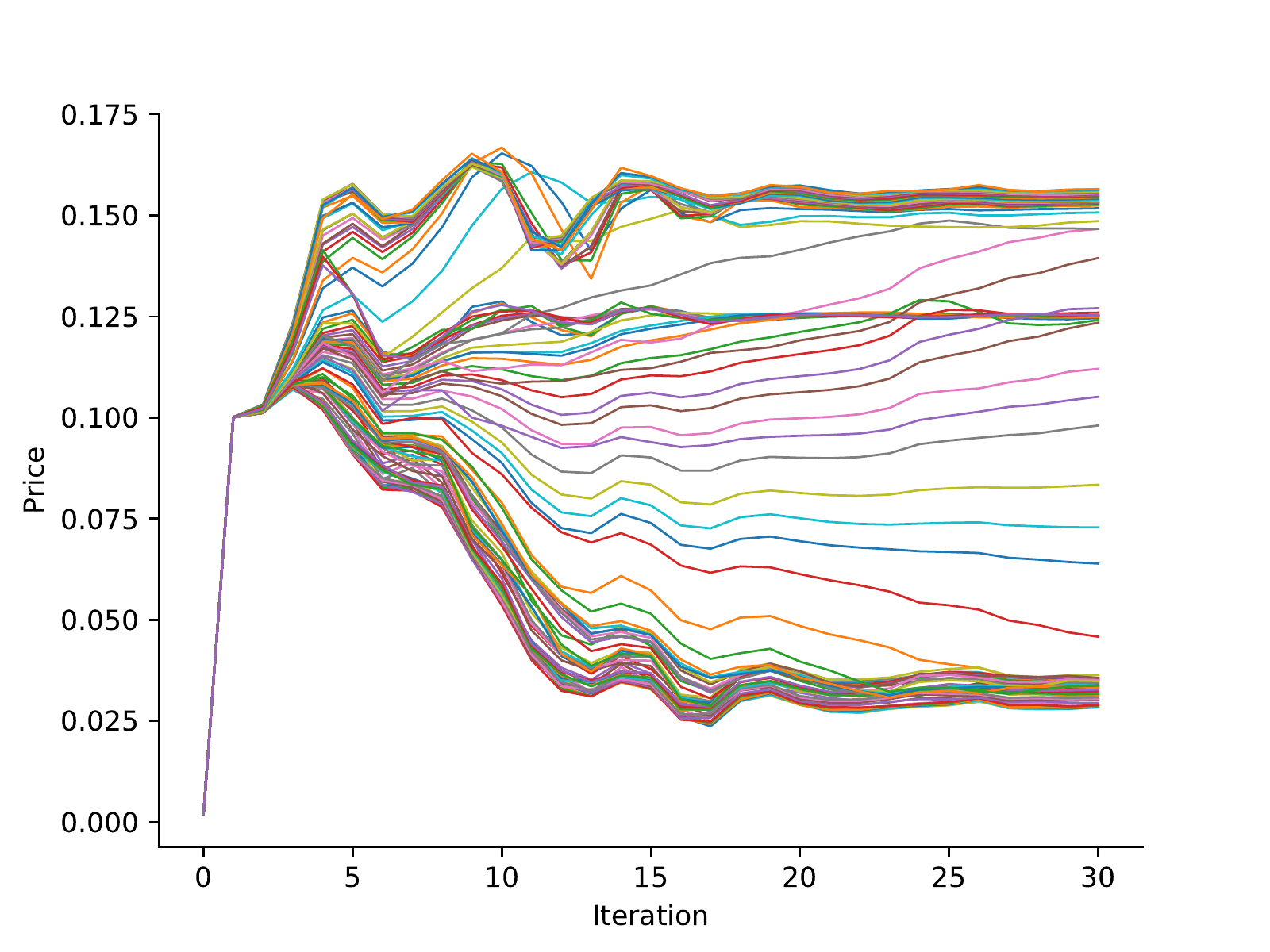} &
    \includegraphics[width=5cm]{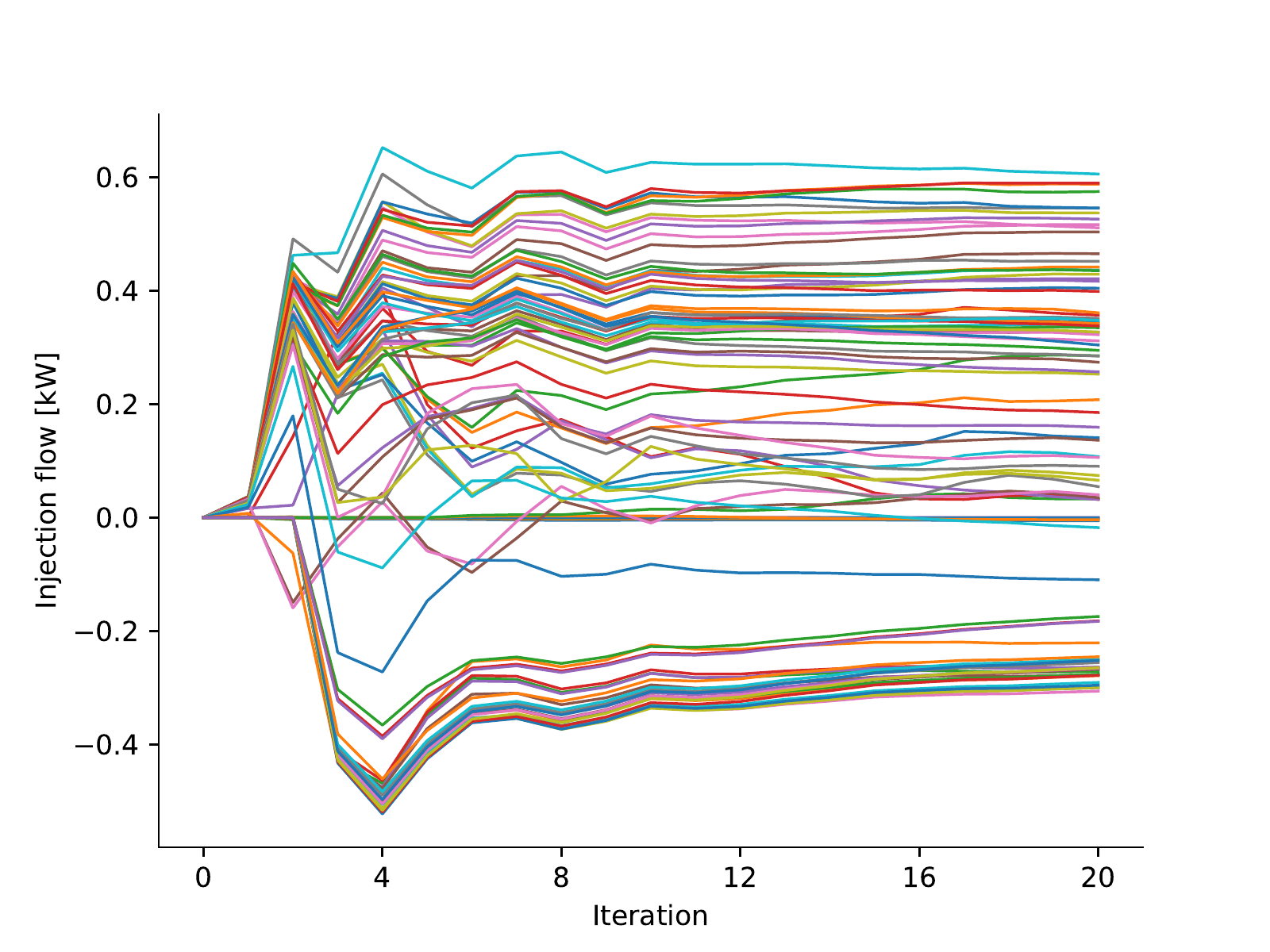} \\
    (a) & (b)
  \end{tabular}
  \caption{Convergence of DADP's prices (a) and PADP resources (b)
           for the 12-Nodes problem}
  \label{fig:numerics:convdadp}
\end{figure}

\paragraph{DADP and PADP convergence.}
We exhibit in Figure~\ref{fig:numerics:convdadp} the convergence
of the DADP's price process and the PADP's resource process along
iterations for the 12-Nodes problem. We depict the convergence only
for the first node, the evolution of price process and resource
process in other nodes being similar.
On the left side of the figure, we plot the evolution
of the 96 different values of the price process
$\price^1 = (\price^1_0, \cdots,  \price^1_{\final-1})$
during the iterations of DADP. We observe that most of the prices
start to stabilize after 15 iterations, and do not exhibit sensitive
variation after 20 iterations. On the right side of the figure, we
plot the evolution of the 96 different values of the resource process
$\alloc^1 = (\alloc_0^1, \cdots, \alloc_{\final-1}^1)$ during the
iterations of PADP. We observe that the convergence of resources
is quicker than for prices, as the evolution of most resources
starts to stabilize after only 10 iterations.

\paragraph{Quality of the exact bounds.}

We then give the lower and upper bounds obtained by SDDP,
DADP, PADP in Table~\ref{tab:district:numeric:upperlower}.
The lower bound of the SDDP algorithm is the value
$\underline{V}_0^{sddp}(x_0)$ given by the SDDP method.
We recall that SDDP returns a lower bound because it uses
a suitable resampling of the global uncertainty distribution
instead of the original distribution itself (see the discussion
in~\S\ref{SDDP_on_the_global_problem}). DADP and PADP lower
and upper bounds are given by Equation~\eqref{eq:district:relaxedmpts}
and Equation~\eqref{eq:district:quantdecompositionrelaxed}
respectively. In Table~\ref{tab:district:numeric:upperlower},
we observe that
\begin{itemize}
  \item SDDP and DADP lower bounds are close to each other,
  \item for problems with more than 12 nodes, DADP's lower
  bound is up to 2.6\% better than SDDP's lower bound,
  \item the gap between the upper bound given by PADP and
  the two lower bounds is rather large.
\end{itemize}
\begin{table}[!ht]
  \centering
  {\normalsize
  \begin{tabular}{|l|ccccc|}
    \hline
    Problem        & \textrm{3-Nodes} & \textrm{6-Nodes} & \textrm{12-Nodes} & \textrm{24-Nodes} & \textrm{48-Nodes} \\
    \hline
    \hline
    SDDP LB        & 225.2            & 455.9            & 889.7             & 1752.8            & 3310.3            \\
    \hline
    DADP LB        & 213.7            & 447.3            & 896.7             & 1787.0            & 3396.4            \\
    \hline
    PADP UB        & 252.1            & 528.5            & 1052.3            & 2100.7            & 4016.6            \\
    \hline
  \end{tabular}
  }
  \caption{Upper and lower bounds given by SDDP, DADP and PADP}
  \label{tab:district:numeric:upperlower}
\end{table}

To sum up, the important result of this paragraph is that,
for optimization problems of large microgrids, DADP is able
to compute a slightly better lower bound than SDDP, and compute
it much faster than SDDP. A parallel version of DADP would
obtain even better performance.

\subsubsection{Policy simulation results}

We now compare the performances of the different algorithms
in simulation. As explained in \S\ref{ssec:nodalpolicies},
we are able to devise online strategies induced by SDDP, DADP
and PADP for the global problem, and to compute by Monte Carlo
an approximation of the expected cost of each of these strategies.

The results obtained in simulation are given
in Table~\ref{tab:district:numeric:simulation}.
SDDP, DADP and PADP values are obtained by simulating the
corresponding strategies on $5,000$ scenarios. The notation
$\pm$ corresponds to the 95\% confidence interval. We use
the value obtained by the SDDP strategy as a reference,
a positive gap meaning that the associated decomposition-based
strategy is better than the SDDP strategy.
Note that all these values correspond to admissible strategies
for the global problem~\eqref{transportproblem}, and thus are
\emph{statistical} upper bounds of the optimal cost~$V_0\opt$
of Problem~\eqref{transportproblem}.

\begin{table}[H]
  \centering
  \resizebox{\textwidth}{!}{
    \begin{tabular}{|l|ccccc|}
    \hline
    Network      & \textrm{3-Nodes} & \textrm{6-Nodes} & \textrm{12-Nodes} & \textrm{24-Nodes} & \textrm{48-Nodes} \\
    \hline
    \hline
    SDDP value   & 226 $\pm$ 0.6    & 471 $\pm$ 0.8    & 936 $\pm$ 1.1     & 1859 $\pm$ 1.6    & 3550 $\pm$ 2.3    \\
    \hline
    \hline
    DADP value   & 228 $\pm$ 0.6    & 464 $\pm$ 0.8    & 923 $\pm$ 1.2     & 1839 $\pm$ 1.6    & 3490 $\pm$ 2.3    \\
    Gap          & - 0.8 \%         & + 1.5 \%         & +1.4\%            & +1.1\%            & +1.7\%    \\
    \hline
    \hline
    PADP value   & 229 $\pm$ 0.6    & 471 $\pm$ 0.8    & 931 $\pm$ 1.1     & 1856 $\pm$ 1.6    & 3508  $\pm$ 2.2   \\
    Gap          & -1.3\%           & 0.0\%            & +0.5\%            & +0.2\%            & +1.2\%    \\
    \hline
  \end{tabular}
  }
  \caption{Simulation results for strategies induced by
           SDDP, DADP and PADP}
  \label{tab:district:numeric:simulation}
\end{table}

We make the following observations.
\begin{itemize}
  \item For problems with more than 6 nodes,
    both the DADP strategy and the PADP strategy beat
    the SDDP strategy.
  \item
    The DADP strategy gives better results than the PADP strategy.
  \item
    Comparing with the last line of Table
    \ref{tab:district:numeric:upperlower}, the statistical
    upper bounds obtained by the three simulation strategies
    are much closer to SDDP and DADP lower bounds than PADP's
    exact upper bound. By assuming that the resource coordination
    process is deterministic in PADP, we impose constant importation
    flows for every possible realization of the uncertainties,
    thus penalizing heavily the PADP algorithm (see also the
    interpretation of PADP in the case of a decentralized
    information structure in \cite[\S3.3]{bounds2019theory}).
\end{itemize}


\section{Conclusion}

In this article, as an application of the companion
paper~\cite{bounds2019theory}, we have studied optimization
problems where coupling constraints correspond to interaction
exchanges on a graph and we have presented a way to decompose
them spatially (Sect.~\ref{chap:district:methods}).
We have outlined two decomposition algorithms, the first relying
on price decomposition and the second on resource decomposition;
they work in a decentralized manner and are fully parallelizable.
Then we have used these algorithms on a specific case study
(Sect.~\ref{chap:district:numerics}), namely the management of
several district microgrids with different prosumers exchanging
energy altogether. Numerical results have showed the effectiveness
of the approach: the price decomposition algorithm beats
the reference SDDP algorithm for large-scale problems with
more than 12~nodes, both in terms of exact bound and induced
online strategy, and in terms of computation time. On problems
with up to 48~nodes (corresponding to 64~state variables), we have
observed that their performance scales well as the number of nodes
grew: SDDP is affected by the well-known curse of dimensionality,
whereas decomposition-based methods are not. Moreover, we have presented
in this article a serial version of the decomposition algorithms,
and we believe that leveraging their parallel nature could decrease
further their computation time.

A natural extension is the following. In this paper, we have only
considered deterministic price and resource coordination processes.
Using larger search sets for the coordination variables, e.g.
considering \emph{Markovian} coordination processes, would make
it possible to improve the performance of the algorithms.
However, one would need to analyze how to obtain a good trade-off
between accuracy and numerical performance.



\bibliographystyle{spmpsci_unsrt}


\end{document}